\documentclass[a4paper,notitlepage, 8pt,reqno]{article}

\usepackage[cp1251]{inputenc}
 \usepackage[T2A]{fontenc}
 \usepackage[russian,english]{babel}
\usepackage[all,arc,poly,2cell,curve,arrow,tips]{xy}
\usepackage{enumerate, eucal, amsthm,amsmath, amssymb}
\usepackage{graphicx}
\usepackage{mathtext}

\usepackage{tikz}

\allowdisplaybreaks[4]

\theoremstyle{definition}
\newtheorem{defn}{Definition}

\theoremstyle{plain}
\newtheorem{thm}{Theorem}
\newtheorem*{thm*}{Theorem}
\newtheorem{prop}{Proposition}
\newtheorem*{prop*}{Proposition}
\newtheorem{cor}{Corollary}
\newtheorem{lem}{Lemma}
\newtheorem{remark}{Remark}
\newtheorem{example}{Example}

\renewcommand{\abstractname}{}

\title{Calculation of  $6j$-symbols for the Lie algebra  $\mathfrak{gl}_n$}

\author{D.V. Artamonov}
\date{}

\begin{document}
    \maketitle

 \maketitle

\renewcommand{\abstractname}{}

\begin{abstract}
	
An explicit description of the multiplicity space that describes  occurrences of  irreducible representations in a splitting of a tensor product of two irreducible finite-dimensional representations of   $\mathfrak{gl}_n$ is given.  Using this description an explicit formula for an arbitrary  $6j$-symbol for finite-dimensional representations  the algebra  $\mathfrak{gl}_n$ is derived. The $6j$-symbol is expressed through a value of a generalized hypergeometric function.

\end{abstract}

In the representation theory of simple Lie algebras there exist natural questions
  concerning a splitting of a tensor product $V\otimes W$ of two irreducible finite dimensional  representations. For example:

\begin{enumerate}
\item	Which irreducible summands $U$ occur in the splitting of   $V\otimes W$ into a direct sum of irreducible finite dimensional  representations? What is a multiplicity of   $U$ (the multiplicity problem)?

\item  What are the explicit formulas for matrix elements of projectors $V\otimes W\rightarrow U$ onto irreducible summands (the problem of calculation of the Clebsh-Gordan coefficients or $3j$-symbols)?

\item What are the matrix elements of an associator, which is an isomorphism between two splittings into irreducible summands of a triple tensor product: $V\otimes (W\otimes U)$ and  $(V\otimes W)\otimes U$ (the problem of calculation of the  Racah coefficients or  $6j$-symbols)?
\end{enumerate}

The problem $1$ should be considered as solved using, for example, the character theory, see the review  \cite{icm}, and also \cite{tt}, \cite{mz}, \cite{kirech}. But there exists a stronger version of the problem  1: the problem of construction of a base in the multiplicity space.

For a long time it was considered that the problems $2$ and $3$ in the general case (i.e. for an arbitrary choice of irreducible finite dimensional representations $V$, $W$, $U$ of  $\mathfrak{gl}_n$, $n\geq 3$)  have no good solutions. But there was a hope that it is possible to obtain a good answer in some  cases if one uses some special functions to express the answer.

Using this idea in  \cite{a1}, \cite{a2} in the case  $\mathfrak{gl}_{3}$  solutions of problems 2 and 3 for an arbitrary choice  of irreducible finite dimensional  $V$, $W$, $U$ were obtained. The answer is not very cumbersome (especially for  $6j$-symbols!) due to the use of  generalized   hypergeometric function.

The key ideas in  the derivation of these results  in  \cite{a1}, \cite{a2} are the following. Firstly, an explicit solution of the multiplicity problem for   $\mathfrak{gl}_3$ that was obtained in \cite{a3} is used.
Secondly, the  so called A-GKZ realization of representations, which is very useful in calculations, is constructed. In the case of   $\mathfrak{gl}_3$ a description of this realization  can be found in    \cite{a1}.

In the present paper we generalize the results of \cite{a2} to the case  $\mathfrak{gl}_n$.  Essentially the scheme  for the calculation of $6j$-symbols from   \cite{a2} is valid in the case  $\mathfrak{gl}_n$ also. But one needs first to solve explicitly the multiplicity problem\footnote{To solve explicitly the multiplicity problem means to construct explicitly a base in a multiplicity space.} and to construct the A-GKZ model for  $\mathfrak{gl}_n$.

The A-GKZ realization of representations of $\mathfrak{gl}_n$ was constructed in \cite{a4}.  So it remains to solve the multiplicity problem. In the present paper we firstly do it in some weak sense  and then we calculate the  $6j$-symbols by analogy with  \cite{a2}.

One should note that there are no papers devoted to calculation of  $6j$-symbols for $\mathfrak{gl}_n$ in the general case  \cite{slp}. Usually some certain cases are considered (see \cite{r82}, \cite{r83}, \cite{r84}, \cite{r85}; note also papers \cite{bl1}, \cite{bl2}, where some classes of  $6j$-symbols are calculated and these coefficients play an important role in the calculation of some Clebsh-Gordan coefficients for  $\mathfrak{gl}_3$).

The problems  1-3  can be posed also  for other series of simple Lie algebras. 
 The problem 1 is considered using the character theory \cite{stf},   \cite{koik}, using the Young tableaux \cite{k}, \cite{g1}, \cite{g2}.
The problems  2 and  3 of calculation of $3j$ and $6j$ symbols are considered only for special cases. Thus in  \cite{al1}, \cite{kl} the $3j$-symbols for symmetric powers of the standard representations were considered (for such representations in the tensor products there are no multiplicities). At the same time the  $6j$-symbols were considered more intensive.  Mostly the
$6j$-symbols for symmetric powers of the standard representations were considered, see  \cite{al2}, \cite{al3}, \cite{al4}, \cite{ju1}, \cite{ho}.
There are also papers where the simplest cases in which non-trivial multiplicities occur are considered, see  \cite{cv}, \cite{cer}, \cite{ju2}. More general cases, as far as I know, were not considered.

There exists a weaker version of the problems  2 and  3:  the problem of algorithmic calculation of  $3j$  and  $6j$-symbols.  This problem is solved completely, see \cite{alg}.

The plan of the present paper is the following.  In Section \ref{predsved} the functional realization of a representation is defined. In this realization using an operation of an overlay of a Young diagram  a new viewpoint to the Weyl construction of an irreducible finite dimensional representation is given. Also a definition of a   $3j$ and a   $6j$-symbol is given.

In Section \ref{if} an explicit solution of the multiplicity problem for a splitting of a tensor product of two irreducible finite dimensional  in the functional realization is given. This problem is equivalent to the problem of description of basic semi-invariants in the triple tensor product.  In such a  form the problem is considered in Section \ref{if}. The main result is the Theorem \ref{tinv}.

In Section  \ref{agkzreal} the basic ideas of the A-GKZ realization  from the paper  \cite{a4} are explained.

In Section \ref{vich6j}  an explicit calculation of a $6j$-symbol is realized. The result is given in Theorem  \ref{ot26j}. Examples of calculations of $6j$-symbols are given.

\section{Preliminary facts.}

\label{predsved}

\subsection{A functional realization}
In the paper Lie algebras and Lie groups over $\mathbb{C}$ are considered. Also we consider only finite dimensional irreducible representations.

Functions on the group $GL_n$ form a representation of the group $GL_n$. Onto a function $f(g)$, $g\in GL_n$, an element  $X\in GL_n$ acts by right shifts

\begin{equation}
\label{xf} (Xf)(g)=f(gX).
\end{equation}

Passing to an infinitesimal version of this action one obtains that the exists an action of $\mathfrak{gl}_n$ in the space of functions on  $G$.

Every irreducible finite dimensional representation can be realized as a sub-representations on the space of functions.
Let $[m_{1},...,m_{n}]$ be a highest weight. In the present paper we suppose that the highest weight is integer and non-negative. This is not an essential restriction, it is done to simplify considerations.

In the space of all functions there exist a highest vector which is written in the following manner.
Let $a_{i}^{j}$ be a function of a matrix element, here $i,j$ run through the sets  of column and row indexes for the group $GL_n$ ($j$ is a row index and $i$ is a column index).
Also put

\begin{equation}
\label{dete}
a_{i_1,...,i_k}:=det(a_i^j)_{i=i_1,...,i_k}^{j=1,...,k},
\end{equation}

where one takes a determinant of a submatrix  in   $(a_i^j)$,
formed by the first rows   $1,...,k$ and columns
$i_1,...,i_k$.

An operator $E_{i,j}$  acts onto a determinant by an action onto column indexes

\begin{equation}
\label{edet1}
E_{i,j}a_{i_1,...,i_k}=\begin{cases}a_{\{i_1,...,i_k\}\mid_{j\mapsto i}}, j\in\{i_1,...,i_k\}\\0 \text{ otherwise }\end{cases},
\end{equation}

where  $.\mid_{j\mapsto i}$ is an operation of a substitution  of  $j$ instead of
$i$.

Take  an integer highest weight  $[m_{1},...,m_{n}]$.
Using  \eqref{edet1}, one can show that the vector
\begin{align}
\begin{split}
\label{stv}
& v_0=a_{1}^{m_{1}-m_{1}}a_{1,2}^{m_{2}-m_{3}}...a_{1,...,n}^{m_{n}}
\end{split}
\end{align}

is a highest vector  for the algebra  $\mathfrak{gl}_n$ with the weight
$[m_{1},...,m_{n}]$.

\begin{thm}[  \cite{zh}]
	\label{lzh}
The space of functions that form a representation with the highest vector 	  \eqref{stv},  is the space of functions that can be written as a polynomial in determinants   $a_{i_1,...,i_k}$, such that their homogeneous powers in determinants of size  $k$ are the same as in the highest vector  \eqref{stv} (i.e. $m_k-m_{k+1}$).
	
\end{thm}

\subsection{ An overlay of a Young symmetrizer}

Let us relate with a highest weight   $[m_{1},...,m_{n}]$ a Young tableau.
It's first row has length $m_{1}$  is filled by   $"1"$, it's second row has length   $m_{2}$ and it is filled by  $"2"$ and so on.  The last row has length $m_{n}$ and it is filled by  $"n"$.
One can relate with this Young tableau a Young symmetrizer, which first performs  an antisymmetrization by columns and then a symmetrization by rows.

The following Proposition takes place that is a direct consequence of the Theorem \ref{lzh}.

\begin{prop}
\label{prp}
If a  monomial in  $a_i^j$ belongs to a representation with the highest vector \eqref{stv} then in every monomial $"1"$ occurs in the set of it's upper indices $m_{1}$ times,
   $"2"$  occurs in the set of it's upper indices  $m_{2}$ times and so on.

\end{prop}
\begin{defn}
{\it An overlay of a Young symmetrizer} onto a monomial in  $a_i^j$, that satisfies the condition of Proposition \ref{prp}, is a result of an application  onto the upper indexes of the monomial of the Young symmetrizer that corresponds to the Young tableau constructed from the highest weight
\end{defn}

\begin{example}

The result of an overlay of a Young symmetrizer onto the monomial  $a_1^1a_2^1a_3^2$ is the following.
$$
a_1^1a_2^1a_3^2+a_1^1a_2^1a_3^2-a_1^2a_2^1a_3^1-a_1^1a_2^2a_3^1=a_1a_{2,3}+a_2a_{1,3}
$$
\end{example}
The following statement can be proved by direct calculations

\begin{prop}
\label{prp1}

Let us be given a monomial in  $a_i^j$, that satisfies the condition of Proposition  \ref{prp}. Then as a result of an overlay of a  Young symmetrizer one obtains a polynomial that belongs to a representation that is described in Theorem \ref{lzh}.
	
Monomials in determinants that satisfy the conditions of Theorem  \ref{lzh} are eigenvectors for an overlay of Young symmetrizer.

\end{prop}

	\subsection{The multiplicity problem}

Take a splitting of a  tensor product of representations  $V$ and $W$  of the algebra $\mathfrak{gl}_n$ into a direct sum of irreducible representations:
\begin{equation}
\label{rzl0}
V\otimes W=\sum_{U} Mult_U\otimes U,    \end{equation}

where    $U$ denotes possible types of irreducible representations that occur in the splitting and $Mult_U$ is the multiplicity space that is a linear space without action of  $\mathfrak{gl}_n$.  One can choose a base  $\{e_f\}$ in this space and put   $U^f:=e_f\otimes U$. Then one can write   \begin{equation}\label{rzl}V\otimes W=\sum_{U,f} U^f.\end{equation}

The multiplicity problem is a problem of construction of a base in the space  $Mult_U$.


\subsection{ Clebsh-Gordan coefficients, $3j$-symbols}

\subsubsection{Clebsh-Gordan coefficients}

Chose in the representations  $V,W,U$  in \eqref{rzl} some bases $\{v_{\alpha}\}$, $\{w_{\beta}\}$, $\{u_{\gamma}\}$. Denote as  $\{u_{\gamma}^f\}$ the corresponding  base  in $U^f$. The Clebsh-Gordan coefficients are numeric coefficients  $D^{U,\gamma,f}_{V,W;\alpha,\beta}\in\mathbb{C}$, that occur in the decomposition

\begin{equation}
\label{kg1}
u_{\gamma}^f =\sum_{\alpha,\beta} D^{U,\gamma,f}_{V,W;\alpha,\beta}  v_{\alpha}\otimes w_{\beta}.
\end{equation}

\subsubsection{$3j$-symbols}

Let us be given representations   $V$, $W$, $U$  of the algebra  $\mathfrak{gl}_n$. Chose in them the bases   $\{v_{\alpha}\}$, $\{w_{\beta}\}$, $\{u_{\gamma}\}$.  A   $3j$-symbol is a collection of numbers

\begin{equation}
\label{3j}
\begin{pmatrix}
V& W& U\\ v_{\alpha} & w_{\beta} & u{\gamma}
\end{pmatrix}^f\in\mathbb{C},
\end{equation}

such that the value

$$
\sum_{\alpha,\beta,\gamma}\begin{pmatrix}
V& W& U\\ v_{\alpha} & w_{\beta} & u{\gamma}
\end{pmatrix}^f    v_{\alpha} \otimes  w_{\beta}  \otimes  u{\gamma}.
$$

is a   $\mathfrak{gl}_n$-semi-invariant. That is this expression is an eigenvector for the Cartan elements  $E_{i,i}$ and vanishes under the action of root elements $E_{i,j}$, $i\neq j$.  The  $3j$-symbols with the same inner indexes form a linear space.  The index   $f$  is numerating basic  $3j$-symbols with the same inner indexes.  One can identify the index  $f$ with a semiivariant that is expressed through the considered   $3j$-symbol.

\subsubsection{A relation between the Clebsh-Gordan coefficients and the $3j$-symbols}
\label{3jsec}

By multiplying   \eqref{rzl} onto a representations  $\bar{U}$, which is contragradient to  $U$ and considering in  $\bar{U}$ a base $\bar{u}_{\gamma}$ dual to  $u_{\gamma}$  one gets a relations

\begin{equation}
\label{3jcg}
D^{U,\gamma,f}_{V,W;\alpha,\beta}=\begin{pmatrix}
V& W& \bar{U}\\ v_{\alpha} & w_{\beta} & \bar{u}_{\gamma}
\end{pmatrix}^f
\end{equation}

Thus the problems of calculation of the Clebsh-Gordan coefficients and the  $3j$-symbols are essentially equivalent.

Also this formula allows to identify the multiplicity spaces for the Clebsh-Gordan coefficients and for the  $3j$-symbols.






 \subsection{The Racah coefficients, $6j$-symbols}

\subsubsection{The Racah coefficients}

The third fundamental problem in the study of tensor products of irreducible representations is the problem of calculation of the Racah coefficients. The Racah coefficients are the matrix elements of the operator that is an isomorphism of  $V^1\otimes (V^2\otimes V^3)$ and $(V^1\otimes V^2)\otimes V^3$. Let us explain this isomorphism in details. A triple tensor product can splitted into a sum of irreducible representations in two ways.

{\bf 1.}  The first way.   Firstly one splits  $V^1\otimes V^2$:

\begin{equation}
\label{r1}
V^1\otimes V^2=\bigoplus_U Mult_{U}^{V^1,V^2}\otimes U,
\end{equation}
where  $U$ is an irreducible representation and  $Mult_{U}^{V^1,V^2}$ is the multiplicity space. Secondly one  multiplies  \eqref{r1} by $V^3$ from the right, and one gets

\begin{equation}
\label{rr1}
(V^1\otimes V^2)\otimes V^3=\bigoplus_{U,W} Mult_{U}^{V^1,V^2}\otimes Mult_{W}^{U,V^3}\otimes W
\end{equation}

{\bf 2.} The second way. Firstly one splits $V^2\otimes V^3$:

\begin{equation}
\label{r2}
V^2\otimes V^3=\bigoplus_U Mult_{H}^{V^2,V^3}\otimes H,
\end{equation}

and then secondly one writes

\begin{equation}
\label{rr2}
V^1\otimes (V^2\otimes V^3)=\bigoplus_{U,W} Mult_{H}^{V^2,V^3}\otimes Mult_{W}^{V^1,H}\otimes W
\end{equation}

There exists an isomorphism  $\Phi: (V^1\otimes V^2)\otimes V^3 \rightarrow V^1\otimes (V^2\otimes V^3)$, which gives a mapping
\begin{equation}
\label{ph}
\Phi: \bigoplus_{U} Mult_{U}^{V^1,V^2}\otimes Mult_{W}^{U,V^3}\rightarrow  \bigoplus_{H} Mult_{H}^{V^2,V^3}\otimes  Mult_{W}^{V^1,H}
\end{equation}

\begin{defn}
The Racah mapping is a mapping  $\Phi$
	
	\begin{equation}
	\label{wrm}
	W\begin{Bmatrix}
	V^1 & V^2 & U\\ V^3 & W &H
	\end{Bmatrix}: Mult_{U}^{V^1,V^2}\otimes Mult_{W}^{U,V^3}\rightarrow  Mult_{H}^{V^2,V^3}\otimes Mult_{W}^{V^1,H}
	\end{equation}

\end{defn}

When one chooses bases in the multiplicity spaces one obtains matrix elements of this mapping. They are called {\it the Racah coefficients}. If $f_1,f_2,f_3,f_4$  are indexes of base vectors in  $Mult_{U}^{V^1,V^2}$, $Mult_{W}^{U,V^3}$, $  Mult_{H}^{V^2,V^3}$,  $Mult_{W}^{V^1,H}$, then one obtains the following notation for the Racah coefficients

\begin{equation}
\label{wrc}
W\begin{Bmatrix}
V^1 & V^2 & U\\ V^3 & W &H
\end{Bmatrix}^{f_1,f_2}_{f_3,f_4}.
\end{equation}

For us it more convenient to deal with close objects called the   $6j$-symbols.
\subsection{$6j$-symbols}

\label{6jsec}


\begin{defn} A $6j$-symbol is a convolution of  $3j$-symbols by the following ruler:
	
	\begin{align}
	\begin{split}
	\label{6js}
	& \begin{Bmatrix}
	V^1 & V^2 & U\\ V^3 & W &H
	\end{Bmatrix}^{f_1,f_2}_{f_3,f_4}:=\sum_{\alpha_1,...,\alpha_6}
	\begin{pmatrix} \bar{V}^1  & \bar{V}^2  & U \\ \bar{v}^1_{\alpha_1} & \bar{v}^2_{\alpha_2}  & u_{\alpha_4}\end{pmatrix}^{f_1}  \cdot \begin{pmatrix}  \bar{U}& \bar{V}^3  &W  \\\bar{ u}_{\alpha_4}  & \bar{v}^3_{\alpha_3} &w_{\alpha_5}   \end{pmatrix}^{f_2}\cdot \\ & \cdot  \begin{pmatrix}  V^2 &  V^3  &  \bar{H} \\ v^2_{\alpha_2} & v^3_{\alpha_3}   & \bar{h}_{\alpha_6}  \end{pmatrix}^{{f}_3}  \cdot \begin{pmatrix} V^1   & H  &\bar{W}  \\ v^1_{\alpha_1}   & h_{\alpha_6} &\bar{w}_{\alpha_5}\end{pmatrix}^{{f}_4}.
	\end{split}
	\end{align}
\end{defn}
Here $\alpha_i$ is an index numerating the base vectors in the corresponding representation.

This expression should be understood as follows: onto a  $3j$-symbol the Lie algebra  $\mathfrak{gl}_n$ acts by acting onto lower indexes. One forms a semi-invariant from four $3j$-symbols using a  convolution of   indexes in such a way that for two  $3j$-symbols only one pair of lower indexes is  convoluted.

There exists the following relation between the Racah coefficients and the $6j$-symbols.
Let us use the fact that there exists a duality between the spaces  $Mult^{V^1,V^2}_U$ and  $Mult^{\bar{V}^1,\bar{V}^2}_{\bar{U}}$. If  $f_1$  is an index of a base vector in  $Mult^{V^1,V^2}_U$ then $\bar{f}_1$ is an index of a dual base vector in $Mult^{\bar{V}^1,\bar{V}^2}_{\bar{U}}$.  Thus  one has

$$
W\begin{Bmatrix}
V^1 & V^2 & U\\ V^3 & W &H
\end{Bmatrix}^{\bar{f}_1,\bar{f}_2}_{f_3,f_4}=\begin{Bmatrix}
V^1 & V^2 & U\\ V^3 & W &H
\end{Bmatrix}^{f_1,f_2}_{f_3,f_4}
$$

In the present paper   below we deal with the  $6j$-symbols only.





\section{The multiplicity problem for the   $3j$-symbols}
\label{if}

In this Section the multiplicity problem for the   $3j$-symbols is solved in the functional realization of representations.
The main result is the Theorem \ref{tinv} which describes functions  $f\in V\otimes W\otimes U$ that are indexing the  $3j$-symbols with the same inner indexes. This theorem is a generalization of an analogous theorem   obtained in \cite{a1} in the case  $\mathfrak{gl}_3$.

In contrast to the case    $\mathfrak{gl}_3$ we do not manage to construct a set independent generators in the space of such functions. The Theorem \ref{tinv} gives a set of linear generators in the space of semi-invariants of a triple tensor product. These generators are indexing   $3j$-symbols with the same inner indexes.

\subsection{ Semi-invariants in $V\otimes W\otimes U$
}

Let us give a description of semi-invariants in  $V\otimes W\otimes U$ in  the functional realization.
Then $V\otimes W\otimes U$  is realized in the space of homogeneous polynomials in matrix elements  $a_i^j, b_i^j, c_i^j$ (the homogeneity conditions are written below in the Proposition   \ref{prp3})   on $GL_n\times GL_n\times GL_n$. Each of these matrix elements can be considered as a vector in the standard vector representation  $V_0\simeq \mathbb{C}^n$ (the algebra  $\mathfrak{gl}_n$ acts onto lower indexes).  Such a viewpoint gives an embedding  $V\otimes W\otimes U\subset V_0^{\otimes T}$, where  $T$ is large enough.

An explicit description of semi-invariants in $V_0^{\otimes T}$ is given essentially by the first principal theorem of the invariant theory  \cite{wey} for the group of lower-unitriangular matrices.

One can consider matrix elements  $a_i^j, b_i^j, c_i^j$ as vectors in different (for different upper indexes, for different symbols $a,b,c$) copies of  $V_0$.  Thus in our notations the first principal theorem can be reformulated as follows.

\begin{thm}
	 \label{ootti}
A semi-invariant for the action of  $\mathfrak{gl}_n$ in the space of polynomials in matrix elements $a_i^j, b_i^j, c_i^j$  is a polynomial in basic semi-invariants that are written as determinants of type   $det(x^{j_1}...x^{j_n})$, where one takes a matrix composed of matrix elements   $x_{i}^j$ where  $x$ is one of the symbols   $a,b,c$ (maybe different for different $j$), the upper index  $j$ takes values $j_1,...,j_n$ and the lower index   $i$ takes values  $1,...,n$.
\end{thm}

 Not all the semi-invariant described in Theorem  \ref{ootti} belong to the functional realization of  $V\otimes W\otimes U$. The following necessary condition, which is a direct corollary of Proposition  \ref{prp}, takes place.

 \begin{prop}
 	\label{prp3}
 Let the highest weights  $V$, $W$, $U$ be $[m_{1},...,m_{n}]$,  $[M_{1},...,M_{n}]$, $[m'_{1},...,m'_{n}]$.
 	
 	If a polynomial in matrix elements $a_i^j, b_i^j, c_i^j$  belongs to   $V\otimes W\otimes U$  then the following condition takes place. Among the variables $a_i^j$ the upper index  $"1"$ occurs $m_{1}$ times,  the upper index   $"2"$ occurs  $m_{2}$  times and so on.
 	Analogous condition must be satisfied by symbols   $b$, $c$.
 \end{prop}

Note that if a polynomial satisfies the conditions of Proposition  \ref{prp3} then one can apply  overlays of the three Young symmetrizes onto  the upper indexes of the symbols  $a$, $b$, $c$.

 The image of the overlays of these three Young symmetrises is the representation $V\otimes W\otimes U$ (see Proposition \ref{prp1})). So if one applies the overlays of these three Young symmetrisers to a semi-invariant that satisfies  the conditions of the Proposition  \ref{prp3}, one gets  a semi-invariant in   $V\otimes W\otimes U$.  Let us describe it explicitly. Introduce an operation of an overlay of an antisymmetrizer   $(\dots)$ onto lower indexes of a monomial in   $a$, $b$, $c$ (also let us call this operation {\it an overlay of brackets} onto lower indexes).

 \begin{defn}
  The operation of an overlay of   $(\dots)$ is defined as follows. One chooses   $n$ lower indexes of a monomial (one says that these indexes are in the brackets) and then one performs an antisymmetrization of these indexes.
 \end{defn}

\begin{example}
	An overlay of  $(\dots)$ onto the first two indexes in the case of the algebra    $\mathfrak{gl}_2$ looks as follows: $a^{1}_{(1}b^{1}_{2)}c^{2}_{1}:= a^{1}_{1}b^{1}_{2}c^{2}_{1}-a_{2}^{1}b^{1}_{1}c^{2}_{1}$
\end{example}

Form the Theorem \ref{ootti}  using an overlay of the Young symmetrizers one obtains the following statement.

\begin{thm}
	\label{tinv}
	Semi-invariants in $V\otimes W\otimes U$ are linear combinations of semi-invariants that are constructed as follows.  One takes a monomial in $a_i^j, b_i^j, c_i^j$, that satisfies thee conditions of Proposition   \ref{prp3}. The lower indexes are divided into non-intersecting groups consisting of    $n$ indexes.
	
	Onto upper indexes the three Young symmetrizes are overlayed. Onto each group consisting of   $n$ lower indexes a bracket  $(\dots)$ is overlayed.
\end{thm}

The following Proposition takes place, which says that the dependence of the constructed semi-invariant on the choice of the overlay of  $(\dots)$ is not so strict.

\begin{prop}
	\label{prp11}
	The function defined in the Theorem \ref{tinv} depends up to sign only on the number of symbols $a$, $b$, $c$ in each bracket  $(\dots)$, but does not depend on the exact placement of upper indexes of these symbols on the first step of the construction of the semi-invariant.
\end{prop}

\proof

Indeed fix a choice of  overlays and let us describe some operations that change these overlays but do not change  up to sign the function. From the existence of these operations the statement of the Proposition follows.

In the formulation below the symbols  $x, y$  denote two different symbols  $a,b,c$, and  $\bullet$  is an arbitrary lower index.

Let us prove the following: one of the anisymmetrizer $(...)$ is overlayed onto a symbol $x^i_{\bullet}$ and  other  anisymmetrizer  $(...)$ is overlayed onto another symbol  $x^j_{\bullet} $, then these symbols can be interchanged. The case $i=j$ is admissible.
	
Indeed suggest first that $i=j$.  Then after application of an overlay of a  Young symmetrizer one obtains an expression that is symmetric in these two indexes.  Thus if one interchanges these symbols  $x^i_{\bullet} $ and then applies the overlay of the antisymetrises one obtains the same expression.

	Now suppose that  $i\neq j$. Without loss of generality one can suggest that these symbols are in one column in the process of  overlay  Young symmetrizers. 
Thus if one interchanges these indexes and then overlays the Young symmetrizers,  one obtains the same expression as one would obtain without interchange of indexes at the beginning but with the sign minus.

\endproof

From the Proposition   \ref{prp}  it follows that for the function constructed on the Theorem   \ref{tinv} one can introduce a notation

  \begin{equation}
 \label{abc1}((x^{j_1}\cdots x^{j_n})\cdots  (y^{i_1}\cdots y^{i_n})),
 \end{equation}

where $x,y,...$  are symbols  $a,b,c$ (the symbols $x$ in $x^{j_1}$, $x^{j_2}$,... can be different symbols $a,b,c$). The upper indexes  must satisfy the conditions of Proposition  \ref{prp3}. The bracket correspond to an overlay of the anitisymmetrization onto lower indexes.

One can also introduce a shorter (but not a complete one) notation
\begin{equation}
\label{abc}
(a^{i_1}\cdots a^{i_{k_1}}b^{j_1}\cdots b^{j_{k_2}}c^{l_1}\cdots c^{l_{k_3}}),
\end{equation}

one must claim that the upper indexes satisfy the conditions of Proposition  \ref{prp3}, and  $k_1+k_2+k_3$ is divisible by  $n$.

\begin{cor}
	A semi-invariant in  $V\otimes W\otimes U$ is a linear combination of semi-invariant of type
	
	\begin{equation}
	\label{3jb}
	f=\prod \frac{1}{t!}(((x^{j_1}\cdots x^{j_n})\cdots  (y^{i_1}\cdots y^{i_n})))^{t},
	\end{equation}
where the exponents  $t$ (which are non-negative integers)  in different factors are different. The set of all upper indexes must satisfy the conditions of Proposition   \ref{prp3}.
\end{cor}

\begin{example}

In the case $\mathfrak{gl}_3$  this construction gives  semiinvariants  defined in\cite{a1}. For example
 $$
(a^1a^2b^1)=det\begin{pmatrix}a_1^1 & a_2^1 & a_3^1 \\ a_1^2 & a_2^2 & a_3^2\\ b_1^1 & b_2^1 & b_3^1 \end{pmatrix},
((c^1c^2b^2)(b^1a^1a^2))=\pm det\begin{pmatrix}a_{2,3} & a_{1,3} & a_{1,2} \\ b_{2,3} & b_{1,3} & b_{1,2}\\ c_{2,3} & c_{1,3} & c_{1,2} \end{pmatrix},
 $$
in  \cite{a1} these semiinvariants are denoted as  $(aab)$ and  $(aabbcc)$.
\end{example}

\begin{example} In the case $\mathfrak{gl}_4$  there exist semiinvariants that are not written in the form of determinants. For example
\begin{align}
\begin{split}
\label{inv4}
&((a^1a^2a^3b^1)(b^2c^1c^2c^3))=a_{1,2,3}b_{2,1}c_{2,3,4}+a_{1,2,3}b_{4,3}c_{4,1,2}-a_{1,2,3}b_{4,2}c_{2,4,1}-\\
&-a_{4,1,2}b_{3,1}c_{2,3,4}+a_{4,1,2}b_{3,4}c_{1,2,3}+a_{4,1,2}b_{3,2}c_{3,4,1}+\\
&+a_{3,4,1}b_{2,1}c_{2,3,4}-a_{3,4,1}b_{2,4}c_{1,2,3}+a_{3,4,1}b_{2,3}c_{4,1,2}.
\end{split}
\end{align}
Since this expression consists of  $9$ terms, it can not be expressed as a determinant.
\end{example}

\section{The A-GKZ realization. The variables $Z$}
\label{agkzreal}
\subsection{The A-GKZ realization of a representation.}

\label{contra}

The A-GKZ realization of a representation of  $\mathfrak{gl}_3$ was introduced in  \cite{a1} for the purpose of calculation of a  $3j$-symbol for this algebra.   A construction of it's analog for   $\mathfrak{gl}_n$ is a non-trivial problem.  It was solved in \cite{a4}. For the purpose of calculation of a  $6j$-symbol it is enough to know only the definition of this realization \footnote{Actually this is just a new viewpoint to the tensor realization of a representation.}. We give it in the present Section.

Consider variables $A_X$ indexed by proper subsets  $X\subset \{1,...n\}$. We claim that  $A_X$ are antisymmetric in    $X$ but they do not satisfy any other relations.   Note that these variables have the same indexes as determinants  \eqref{dete}, but   $A_X$ do not satisfy other relation but antisymmetry.

Onto these variables the algebra acts by the ruler

\begin{equation}
\label{edet10}
E_{i,j}A_{i_1,...,i_k}=\begin{cases}
A_{\{i_1,...,i_k\}\mid_{j\mapsto i}},\text{  if  }j\in \{i_1,...,i_k\}\\ 0 \text{ otherwise }.
\end{cases}
\end{equation}

and onto the product of these variables the algebra acts according to the Leibnitz ruler. Thus the algebra of polynomials    $\mathbb{C}[A]$ is a representation of  $\mathfrak{gl}_n$.

Consider the system of partial differential equations called the A-GKZ system, which is constructed as follows.
Let  $I\subset \mathbb{C}[A]$ be an ideal of relations between the determinants $a_X$. It is known that it is  generated by the Plucker relations.

Under the mapping

$$
A_X\mapsto \frac{\partial}{\partial A_X},
$$
  the ideal  $I$ is transformed to an ideal   $\bar{I}\subset \mathbb{C}[\frac{\partial}{\partial A}]$ in the ring of differential operators with constant coefficients.

The A-GKZ system is a system of partial differential equations defined by the ideal   $\bar{I}$: $$\forall \mathcal{O}\in \bar{I}:\,\,\,\,\mathcal{O}F(A)=0.$$ One has.

\begin{thm} [\cite{a4}]  The space of polynomial solutions of the A-GKZ system is a representation of  $\mathfrak{gl}_n$. This representations is a direct sum of all finite dimensional irreducible representations with an integer highest weight taken with the multiplicity   $1$.
	
	A sub-representation with the highest weight   $[m_1,...,m_n]$ is a space of polynomial solutions such their the  homogeneous power of $A_{X}$ with   $|X|=1$ equals  $m_1-m_2$, homogeneous power of $A_{X}$ with   $|X|=2$  equals $m_2-m_3$ and so on.
\end{thm}

This Theorem gives a realization of finite dimensional irreducible representations, which is called the A-GKZ realization.

Note that the substitution

$$
A_X\mapsto a_X
$$

maps isomorphically the A-GKZ realization to the functional realization form the Theorem  \ref{lzh}.

Define the action

\begin{equation}
\label{dve}
f(A)\curvearrowright  g(A):=f(\frac{d}{dA})g(A),
\end{equation}
then on the space of polynomials in variables   $A_X$ there exists an invariant\footnote{A scalar product is invariant if  $<E_{i,j}f,g>=-<f,E_{j,i}g>$} scalar product

\begin{equation}
\label{skd}
<f(A),g(A)>= f(A)\curvearrowright  g(A)\mid_{A=0}.
\end{equation}

Due to the symmetry of the scalar product one can also write
 $ <f(A),g(A)>= g(A)\curvearrowright  f(A)\mid_{A=0}.$

Note that if the representation $V$ is realized  in the space of polynomials in the variables $A_X$ (i.e. $V$ is some space consisting of polynomials $ \{h(A)\}$), then the contragradient representation is realized in the space of polynomials in the operators $\frac{\partial}{\partial A_X}$.  The action  of $\mathfrak{gl}_n$ on the differential operators is generated by an the action on functions of $A_X$. In this case $\bar{V}=\{h(\frac{\partial }{ \partial A_X}): \,\,\, h(A)\in V\}$. The pairing is given by a formula similar to \eqref{skd}:

\begin{equation}
	\label{spar}
	<h_1(A),h_2(\frac{\partial }{ \partial A_X})>=h_2(\frac{\partial }{ \partial A_X})h_1(A)\mid_{A=0}.
\end{equation}

\subsection{Semi-invariants in the A-GKZ realization}

Let $f(a,b,c)$  be a semi-invariant in   $V\otimes W\otimes U$ constructed in the present paper. In the A-GKZ realization one can consider a polynomial $f(A,B,C)$,  where the symbols  $a,b,c$ are changed to  $A,B,C$. The obtained polynomial in general does not belong to  $V\otimes W\otimes U$ in the A-GKZ realization.



But the functional and the A-GKZ realizations are realizations of the same representation (in our case of this is the triple tensor product of irreducible representations). Suppose that to the vector $f(a,b,c)$ in the A-GKZ representation there corresponds a vector
$F(A,B,C)$ in the functional representation. Since the A-GKZ realization is transformed into the functional realization when one imposes the Plucker relation, the following holds

\begin{equation}
F(A,B,C)=f(A,B,C)+\sum_{\beta} pl^A_{\beta}f^1_{\beta}+ pl^B_{\beta}f^2_{\beta}+ pl^C_{\beta}f^1_{\beta},
\end{equation}
where  $pl^A_{\beta}$  are the basic Plucker relations for the variables  $A_X$ and  $f^1_{\beta}$ is some polynomial in the variables   $A, B, C$.

One gets that
\begin{align}
\begin{split}
\label{gsk}
<F(A,B,C),F_{\mu}(A)F_{\nu}(B)F_{\nu}(C)>=<f(A,B,C),F_{\mu}(A)F_{\nu}(B)F_{\nu}(C)>.
\end{split}
\end{align}

where  $F_{\mu}(A)$, $F_{\nu}(B)$, $F_{\nu}(C)$  are solutions of the A-GKZ system.

\subsection{The variables  $Z$, the numbers $z_{\alpha}$}
\label{razz}


Consider the semi-invariant   \eqref{abc1}. Since we are considering the functional realization   one can write it as a function of determinants
\begin{equation}
\label{slz}
\sum_{\alpha} z_{\alpha}a_{}^{p_{\alpha}}b^{q_{\alpha}}_{}c_{}^{r_{\alpha}},
\end{equation}

where  $\alpha$ is an index numerating the summands in the explicit expression for   \eqref{3jb}, and  $z_{\alpha}$ is a numeric coefficient. One understands  $a_{}^{p_{\alpha}}$, $b^{q_{\alpha}}_{}$, $c_{}^{r_{\alpha}}$  using the multi-index notation, that is   $a_{}^{p_{\alpha}}=\prod_{X}a_{X}^{p_{\alpha,X}}$.

Introduce variables that correspond to summands in the obtains sum. One has a natural notation for them  $Z_{\alpha}=[a_{}^{p_{\alpha}}b^{q_{\alpha}}_{}c_{}^{r_{\alpha}}]$.  The set of the obtained variables (for all possible  \eqref{abc1})  denote as $Z$:

\begin{equation}
\label{pz}
Z=\{ Z_{\alpha}= [a_{}^{p_{\alpha}}b^{q_{\alpha}}_{}c_{}^{r_{\alpha}}],...   \}
\end{equation}

\begin{example}
If one takes  ta function $f$ of type  \eqref{inv4}, then the collection of variables  $Z$ looks as follows

\begin{align*}
&Z=\{[a_{1,2,3}b_{2,1}c_{2,3,4}],[a_{1,2,3}b_{4,3}c_{4,1,2}],[a_{1,2,3}b_{4,2}c_{2,4,1}],[a_{4,1,2}b_{3,1}c_{2,3,4}],[a_{4,1,2}b_{3,4}c_{1,2,3}],\\&[a_{4,1,2}b_{3,2}c_{3,4,1}],[a_{3,4,1}b_{2,1}c_{2,3,4}],[a_{3,4,1}b_{2,4}c_{1,2,3}],[a_{3,4,1}b_{2,3}c_{4,1,2}]\}.
\end{align*}

The coefficients  $z_{\alpha}$ are equal to numbers  $\pm 1$,  occuring at the corresponding summd in  \eqref{inv4}.

\end{example}

Note that the exists a natural mapping

\begin{equation}
\label{zx}
Z_{\alpha}= [a_{}^{p_{\alpha}}b^{q_{\alpha}}_{}c_{}^{r_{\alpha}}]\mapsto z_{\alpha}\in\mathbb{C}
\end{equation}

One can consider  $f$  of type \eqref{3jb} as a polynomial in variables  $Z$.

\begin{defn}
	Define a support of function written as a power series as a set of exponents of the involved monomials. Denote the support as   $suppf$.
\end{defn}

Since $f$ is of type \eqref{3jb}, one has.

\begin{lem}
	\label{lrb} For the support of the function $f$ considers as a function of   $Z$, one has
	\begin{equation}
	\label{rb}
	suppf=(\kappa+B)\cap \text{(the non-negative octant)}
	\end{equation}
 for some constant vector  $\kappa$ and some lattice $B$.
\end{lem}
\proof

Let us give first the construction of the lattice  $B$.  Take a factor of type \eqref{abc1}  in  \eqref{3jb}. When one defines this factor one fixes on overlay of brackets    $(\dots)$ onto lower indexes. In this procedure one substitutes  into the lower indexes the numbers $1,...,n$.   When one fixes a substitution into each bracket one obtains a variable   $Z_{\alpha}$.  With each such a variables one relates  a unit vector $e_{Z_{\alpha}}$ in the space of exponents of monomials in variables  $Z$.  Take vectors  $e_{Z_{\alpha}}-e_{Z_{\beta}}$ for all possible pairs of variables   $Z_{\alpha}$, $Z_{\beta}$ from one factor of type   \eqref{abc1}, for all possible factors of type   \eqref{abc1}.
The lattice  $B$ is generated by these differences.

The vector  $\kappa$  is defined as a vector of exponents of a monomial in variables   $Z$, which appears if one fixes  in
a decomposition \eqref{slz} of  \eqref{abc1}. Such a fixation is done for all factors of type  \eqref{abc1}, occurring in $f$ of type \eqref{3jb}.

By construction  $suppf\subset (\kappa+B)\cap \text{(the non-negative octant)}$.  One needs to prove the coincidence of these sets.

By definition  $b\in B$ is a shift of a vector of exponents when one changes a substitution of $1,...,n$ in a bracket  $(\dots)$ in the construction of \eqref{abc1}.

But  $(\dots)$ is an antisymmetrization over {\it all }  possible substitutions of   $1,...,n$. Thus arbitrary shifts from the initial vector in the case when one obtains a vector with non-negative coordinates are vectors from $suppf$.


\endproof

\begin{remark}
	The lattice  $B$ in  \eqref{rb} is actually defined by the collection of variables  $Z$.  And the function   $f$  defines the initial vector   $\kappa$ in  \eqref{rb}.
\end{remark}

There exists a mapping

\begin{align*}
&[a_{}^{p_{\alpha}}b^{q_{\alpha}}_{}c_{}^{r_{\alpha}}]\mapsto A_{}^{p_{\alpha}}B^{q_{\alpha}}_{}C_{}^{r_{\alpha}}
\end{align*}

from the space of polynomials in variables  $Z$  to the space of polynomials in variables  $A$, $B$, $C$.
Denote as   $pr_A$, $pr_B$, $pr_C$  the induced mappings from the space of exponents of variables   $Z$ to the space of exponents of variables  $A$, $B$, $C$.

\section{The $6j$-symbols}

\label{vich6j}

Now let us calculate an arbitrary  $6j$-symbol for the algebra  $\mathfrak{gl}_n$.  Let us follow the scheme of calculation of a  $6j$-symbol for the algebra   $\mathfrak{gl}_3$  from \cite{a2}.   The considerations on the present paper follow literally the considerations from  \cite{a2} until the the construction of the function  \eqref{fj}. But then in the formulation of the Theorem  \ref{ot26j} there is an important difference from the Theorem 4 in  \cite{a2}. In  \cite{a2} in an explicit expression of a basic semi-invariants as polynomial in determinants the coefficients at monomials are equal to $\pm 1$. Thus in  \cite{a2}, into the functions \eqref{fj}  in the Theorem 4 one substitutes  $\pm 1$.  For the semi-invariants considered in the present paper  these coefficients can take values denoted as   $z_{\alpha}$, that are not necessarily  $\pm 1$. Thus into a function  \eqref{fj}  in the Theorem    \ref{ot26j} one substitutes other values.

\subsection{An expression through the convolution}

\begin{lem}
	\begin{align}
	\begin{split}
	\label{6jss}
	&	\begin{Bmatrix}
	V^1 & V^2 & U\\ V^3 & W &H
	\end{Bmatrix}^{f_1,f_2}_{f_3,f_4}= f_1(  \frac{\partial}{\partial A^1}, \frac{\partial}{\partial A^2} , A^4 )f_2(  \frac{\partial}{\partial A^4}, \frac{\partial}{\partial A^3} , A^5 )\cdot\\&\cdot f_3(  \frac{\partial}{\partial A^2},  A^3, \frac{\partial}{\partial A^6}  )f_4(  A^1,  A^6, \frac{\partial}{\partial A^5}  ).\mid_{A^1=...=A^6=0}
	\end{split}
	\end{align}
\end{lem}

Here  $A^i$, $i=1,...,6$ are six copies of independent sets of variables  $A_X^i$, where  $X\subset \{1,...,n\}$ are proper subsets. As usual  these variables are anisymmetric in $X$ but do not satisfy any other relations.

The proof of this Lemma in the case   $\mathfrak{gl}_n$ is literally the same as in the case  $\mathfrak{gl}_3$  in  \cite{a2}.  But let us write it here.

\proof
Let us  use the formula \eqref{6js}.
We need to calculate  $3j$-symbols for the contragradient representation and the dual basis. Let use a realization of a contragradient representation described at the end of the section \ref{contra}. Suppose that we take an orthogonal base $F_{\alpha_i}(A^i)$. The a base dual to $F_{\alpha_i}(A^i)$ is  $\frac{1}{|F_{\alpha_i}|^2}F_{\alpha_i}(\frac{\partial}{\partial A^i})$.

Note that a $3j$-symbol of the form  $$\begin{pmatrix} \bar{V}^1 & \bar{V}^2 & U\\F_{\alpha_1}(\frac{\partial}{\partial A^1}) &  F_{\alpha_2}(\frac{\partial}{\partial A^2}) & F_{\alpha_4}(A^4)   \end{pmatrix}^f$$
can be calculated as follows:

\begin{equation}
	\label{s3j1}
	\begin{pmatrix} \bar{V}^1 & \bar{V}^2 & U\\F_{\alpha_1}(\frac{\partial}{\partial A^1}) &  F_{\alpha_2}(\frac{\partial}{\partial A^2}) & F_{\alpha_4}(A^4)   \end{pmatrix}^f=\frac{<f(\frac{\partial}{\partial A^1},\frac{\partial}{\partial A^2},A^4),F_{\alpha_1}(\frac{\partial}{\partial A^1}) F_{\alpha_2}(\frac{\partial}{\partial A^2}) F_{\alpha_4}(A^4)>}{|F_{\alpha_1}(\frac{\partial}{\partial A^1})|^2| F_{\alpha_2}(\frac{\partial}{\partial A^2})|^2 |F_{\alpha_4}(A^4)|^2}.
\end{equation}

The scalar product in case when there is a function not of  a variable, but of a differentiation operator is calculated using a formula similar to \eqref{skd}. One has

\begin{align}
	\begin{split}
		\label{ss1}
		&<f(\frac{\partial}{\partial A^1},\frac{\partial}{\partial A^2},A^4),F_{\alpha_1}(\frac{\partial}{\partial A^1}) F_{\alpha_2}(\frac{\partial}{\partial A^2}) F_{\alpha_4}(A^4)>=
		<f( A^1, A^2,A^4),F_{\alpha_1}(A^1) F_{\alpha_2}( A^2) F_{\alpha_4}(A^4)>,\\
		& |F_{\alpha_1}(A^1)|^2=|F_{\alpha_1}(\frac{\partial}{\partial A^1})|^2,...
	\end{split}
\end{align}

Bases $F_{\alpha_1}(A^1)$ and $F_{\alpha_1}(\frac{\partial}{\partial A^1})$ etc. are not dual,  the basis  dual to $F_{\alpha_1}(A^1)$ is  $\frac{1}{|F_{\alpha_1}|^2}F_{\alpha_1}(\frac{\partial}{\partial A^1})$. So the $6j$-symbol is expressed in terms of the considered $3j$-symbols \eqref{s3j1} as follows

\begin{align}
	\begin{split}
		\label{6js1}
		& \begin{Bmatrix}
			V^1 & V^2 & U\\ V^3 & W &H
		\end{Bmatrix}^{f_1,f_2}_{f_3,f_4}:=\sum_{\alpha_1,...,\alpha_6}
		\begin{pmatrix} \bar{V}^1 & \bar{V}^2 & U\\F_{\alpha_1}(\frac{\partial}{\partial A^1}) &  F_{\alpha_2}(\frac{\partial}{\partial A^2}) & F_{\alpha_4}(A^4)   \end{pmatrix}^{f_{1}}  \cdot\\&\cdot \begin{pmatrix}  \bar{U}& \bar{V}^3  &W  \\   F_{\alpha_4}(\frac{\partial}{\partial A^4})  &   F_{\alpha_3}(\frac{\partial}{\partial A^3}) & F_{\alpha_4}(A^5) \end{pmatrix}^{f_2}\cdot \\ & \cdot  \begin{pmatrix}  V^2 &  V^3  &  \bar{H} \\  F_{\alpha_2}(A^2) & F_{\alpha_3}(A^3)   & F_{\alpha_6}(\frac{\partial}{\partial A^6}) \end{pmatrix}^{{f}_{3}}  \cdot \begin{pmatrix} V^1   & H  &\bar{W}  \\  F_{\alpha_1}(A^1)   &  F_{\alpha_6}(A^6)  &  F_{\alpha_5}(\frac{\partial}{\partial A^5})\end{pmatrix}^{{f}_{4}} \cdot |F_{\alpha_1}|^2\cdot....\cdot |F_{\alpha_6}|^2
	\end{split}
\end{align}

Take the expressions \eqref{s3j1} and substitute them in \eqref{6js1}. Consider \eqref{ss1}. At the same time the expression $|F_{\alpha_i}|^2$  occurring
at the end of \eqref{6js}   are written as $F_{\alpha_i}(\frac{\partial}{\partial A^i})F_{\alpha_i}(A^i)\mid_{A=0}$. 

One obtains

 \begin{align*}
	&В   \begin{Bmatrix}
		В  V^1 & V^2 & U\\ V^3 & W &H
		В  \end{Bmatrix}^{f_1,f_2}_{f_3,f_4}=\sum_{\alpha_1,...,\alpha_6}\frac{<f_{1},F_{\alpha_1}F_{\alpha_2}F_{\alpha_4}>}{  |F_{\alpha_1}|^2     |F_{\alpha_2}|^2   |F_{\alpha_4}|^2 }   F_{\alpha_1}(\frac{\partial}{\partial A^1})F_{\alpha_2}(\frac{\partial}{ \partial A^2})F_{\alpha^4}(A^4) \cdot\\& \frac{<f_{2},F_{\alpha_4}F_{\alpha_3}F_{\alpha_5}>}{  |F_{\alpha_4}|^2     |F_{\alpha_3}|^2   |F_{\alpha_5}|^2 }   F_{\alpha_4}(\frac{\partial}{\partial A^4})F_{\alpha_3}(\frac{\partial}{ \partial A^3})F_{\alpha_5}(A^5)...\mid_{A_1=...=A_6=0}
	В  \end{align*}

Now write

$$ f_1(  \frac{\partial}{\partial A^1}, \frac{\partial}{\partial A^2} , A^4 )=\sum  \frac{<f_{1},F_{\alpha_1}F_{\alpha_2}F_{\alpha_4}>}{  |F_{\alpha_1}|^2     |F_{\alpha_2}|^2   |F_{\alpha_4}|^2 }   F_{\alpha_1}(\frac{\partial}{\partial A^1})F_{\alpha_2}(\frac{\partial}{ \partial A^2})F_{\alpha^4}(A^4),
$$
and analogous expressions for$ f_2(  \frac{\partial}{\partial A^4}, \frac{\partial}{\partial A^5} , A^5 )$ , $f_3(  \frac{\partial}{\partial A^2},  A^3, \frac{\partial}{\partial A^6}  )$,   $f_4(  A^1,  A^6, \frac{\partial}{\partial A^5}  )$.

Using that  $\{F_{\alpha_i}(A^i),   F_{\alpha'_i}(\frac{\partial}{\partial A^i}) \}=|F_{\alpha_i}|^2$, if   $\alpha_i=\alpha'_i$ and  $0$ otherwise one gets the statement of the Lemma.

\endproof

\subsection{The selection ruler, the lattice  $D$}

\label{poror}

In the formula \eqref{6jss}  for a $6j$-symbol the functions  $f_1,f_2,f_3,f_4$ are involved.  Let us substitute into each  $f_i$ instead of a differential operator the corresponding variable. Then one can consider the functions $f_i$ in two ways. 

Firstly in Section  \ref{razz} there was introduces a collection of variables  $Z$ and the function  $f$ of type \eqref{3jb} was considered as a functions of these variables.  Consider the functions  $f_i$  as functions of there own collections of variables  $Z^1,Z^2,Z^3,Z^4$.
In this case the support of  $f_i$ belongs to some space   $\mathbb{Z}^M$ (where  $M$ is the number of variables in $Z^i$).

Secondly, one can consider $f_i$ as a function of variables   $A^j$ (the set of indices  $j$for the variables $A$ involved in   $f_i$ is taken form the formula   \eqref{6jss}). We consider as different the variables $A^j$ involved in different  $f_i$, thus let us introduce  a notation $A^j_{X,i}$ for a variable  $A^j_X$, involved in  $f_i$.  According to this approach the support of  $f_i$ belongs to the space $\mathbb{Z}^m$ (here  $m=3(2^n-2)$).  Note that in the space $\bigoplus_{i=1}^4(\mathbb{Z}^{m})$ one can introduce a base vector $e^{A^j_{X,i}}$.


We have defined the projections  $pr^i$ from  $\mathbb{Z}^M$ to $\mathbb{Z}^m$,  induced by natural substitutions of variables   $A^j$ instead of the variables  $Z^i$.
Also let

\begin{equation}\label{pr}
pr:=\oplus_{i=1}^4 \\
pr^i:\bigoplus^4\mathbb{Z}^M\rightarrow \bigoplus^4\mathbb{Z}^m
\end{equation}

We proved in Lemma \ref{lrb} that if one considers  $f_i$ as a function of the variables  $Z^i$,  then for it's support  $supp_{Z^i}f_i\subset \mathbb{Z}^M$ one has

\begin{equation}
\label{spf}
supp_{Z^i}f_i=(\kappa_i+B) \cap  \mathbb{Z}^{M}_{\geq 0}
\end{equation}

Introduce a notation

$$
H:=supp_{Z^1}f_1\oplus  supp_{Z^2}f_2\oplus supp_{Z^3}f_3\oplus supp_{Z^4}f_4,
$$

note that $H$  is an intersection of the non-negative octant and the shifted lattice  
$(\kappa_1\oplus\kappa_2\oplus\kappa_3\oplus \kappa_4)+B\oplus B\oplus B\oplus B$.  Also $H$  is a suuport of the fucntion   $f_1\cdot...\cdot f_4$ as a function of variables $Z^1,...,Z^4$.

Now introduce a lattice   $D$. 

\begin{defn}
	\label{d} 
According to substitution of arguments into $f_i$\footnote{remind that additionally we substitute into each  $f_i$  indtead of a differential operator a variable.} in \eqref{6jss}, define the lattice
	$
	D\subset \bigoplus_{i=1}^4 (\mathbb{Z}^{m})
	$
	as lattice generated for all possible  $X\subset\{1,...,n\}$   by the vectors that are sums of $e^{A^j_{X,i}}$, отвечающих двум координатам corresponding to coordinates with  the same   $X$  and  the  same variables  $A^j_X$ but different    $i$.

\end{defn}

Thus the lattice  $D$ is generated by the vectors

 \begin{align*}
& e^{A^1_{X,1}}+e^{A^1_{X,4}}, \,\,\,     e^{A^2_{X,1}}+e^{A^2_{X,3}},   \,\,\,      e^{A^3_{X,2}}+e^{A^3_{X,4}}, \\& e^{A^4_{X,1}}+e^{A^4_{X,2}} ,   \,\,\,    e^{A^5_{X,2}}+e^{A^5_{X,4}},  \,\,\,   
  e^{A^6_{X,3}}+e^{A^6_{X,4}} 
 \end{align*}

Let us be given a monomial that is the decomposition of  \eqref{6jss}. It gives a non-zero input if it satisfies the following condition. For every variable   $A^j_X$, $j=1,...,6$  the order of differentiation by the variable  $A^j_X$  equals to the exponent of the variable  $A^j_X$. The condition of existence of such monomials reformulated in terms of the supports of functions in variables  $Z^1,...,Z^4$ gives the  following statement.

\begin{thm}[The selection ruler] \label{ot1}
	
If the $6j$-symbol   \eqref{6jss} is non-zero then
	
	$$
	H\cap pr^{-1}( D)\neq \emptyset,
	$$
\end{thm}

\subsection{The formula for a   $6j$-symbol}

Let us proceed to the calculation of  \eqref{6jss}. 
The procedure of calculation is the following.

Consider  $f_1...f_4$ as functions of variables $Z^1,Z^2,Z^3,Z^4$.  Present  $f_1\cdot...\cdot f_4$ as sums of monomials in these variables. Note that we have now coefficients   $z_{\alpha}^i$,  $i=1,...,4$.  Take only the summands whose exponents belong to    $
H\cap pr^{-1}( D)\neq \emptyset
$. Change all the variables from   $Z^1,Z^2,Z^3,Z^4$ to  $A^j_X$  or  $\frac{\partial}{\partial A_X^j}$  according to the arguments of   $f_i$ in \eqref{6jss}.  Multiply them assuming that the variables and the differential operators commute. Then in all the obtained monomials in $A^j_X$  and $\frac{\partial}{\partial A_X^l}$ one applies the differential operators to variables and then substitutes into the variables zero.

Thus if for example one considers the variable  $A_1^1$ (i.e. $X=\{1\}$), then our actions look as follows. Such a symbol occurs  in the variables in the collections   $Z^1$ and  $Z^4$.  Take a monomial that is obtained in the decomposition of   $f_1\cdot...\cdot f_4$.  Let it's support belong to  $
H\cap pr^{-1}( D)$.
Write it explicitly with a coefficient in this monomial. This coefficient is a product of inverse values of factorials of exponents, that come from  \eqref{3jb}, together with a numeric coefficient of type $z_{\alpha}$,  that occurs at the variable form the collection  $Z=\{Z^1,...,Z^4\}$ in the decomposition \eqref{skd} of the factors in  \eqref{3jb}:

\begin{equation}
\label{monom}
\underbrace{\frac{(z^1_{\alpha_1}[A^1_1...])^{\beta_1}}{\beta_1!}     \frac{(z^1_{\alpha_2}[A^1_1...])^{\beta_2}}{\beta_2!}...         }_{\text{from }f_1}\cdot...\cdot
\underbrace{\frac{(z^4_{\delta_1}[A^1_1...])^{\gamma_1}}{\gamma_1!}     \frac{(z^4_{\delta_2}[A^1_1...])^{\gamma_2}}{\gamma_2!}...         }_{\text{from }f_4}.
\end{equation}

Then one calculates the sum of exponents of variables, whose notation contains  $A_1^1$.  For the factors originating from  $f_1$ this sum equals    $\beta_1+\beta_2+..$, and for the factors originating from   $f_4$ this sum equals   $\gamma_1+\gamma_2+...$. The fact that the supports belong to  $
H\cap pr^{-1}( D)$ implies that   $\beta_1+\beta_2+...=\gamma_1+\gamma_2+...$.
When one passes to $A^j_X$  or  $\frac{\partial}{\partial A_X^j}$   into the factors originating from  $f_1$,  one substitutes   $\frac{\partial}{\partial A_1^1}$, and into the factors originating from  $f_4$ one substitutes   $A_1^1$. After the application of the differential operator to   $A_1^1$ and  the substitution of zero into $A_1^1$ one essentially removes from \eqref{monom} all the symbols $A_1^1$ and writes in the top the factor  $(\beta_1+\beta_2+...)!$.

When one performs analogous operations with all the variables  $A_X^j$ the monomial \eqref{monom} is transformed to a numeric fraction. It's denominator is a factorial  (in the multi-index sense) of the exponent of the monomial considered as a monomial in variables   $Z^1,...,Z^4$, and in the numerator is a factorial of the exponent (in the multi-index sense) of the monomial considered as a monomial in variables  $A_X^j$.  The obtained fraction must be multiplied by  $(z^1_{\alpha_1})^{\beta_1}\cdot (z^1_{\alpha_2})^{\beta_2}\cdot...$

Let us give a formula for a  $6j$-symbol.  The set $
H\cap pr^{-1}( D)$ is a shifted lattice in the space of exponents of monomials in variables  $Z^1,...,Z^4$. Hence for some vector   $\varkappa$ and some lattice  $L\subset  (\mathbb{Z}^M)^{\oplus 4}$ one can write

\begin{equation}
\label{kp}
H\cap pr^{-1}( D)=\varkappa+L\subset  (\mathbb{Z}^M)^{\oplus 4}
\end{equation}

There exists a projection  $pr$,  defined by the formula  \eqref{pr}.  Let us related with the shifted lattice    $\varkappa+L$ a hypergeometric type series (which in fact is a finite sum) in the variables  $\mathcal{Z}=\{Z^1,...,Z^4\}$, defined by the formula:

\begin{equation}
\label{fj}
\mathcal{J}_{\gamma}(\mathcal{Z};L)=\sum_{x\in \varkappa+L} \frac{\Gamma(pr(x)+1)\mathcal{Z}^{x}}{\Gamma(x+1)}
\end{equation}

\begin{thm}\label{ot26j}
The	$6j$-symbol \eqref{6jss} equals $\mathcal{J}_{\gamma}(z;L)$, where instead of a variable from the collection    $\mathcal{Z}=\{Z^1,...,Z^4\}$ one substitutes the number  $z_{\alpha}$   by the ruler \eqref{zx}.
	
\end{thm}

 \subsection{Example of Calculation}

Consider the algebra $\mathfrak{gl}_4$. To define the $6j$-symbol, we first fix the semi-invariants:  

$$
f_1=(aabc), \quad f_2=(abbc), \quad f_3=(abbc), \quad f_4=(aabc).  
$$

In expression \eqref{6jss} for the $6j$-symbol, we must substitute the variables $A^j_X$ or the operators  $\frac{\partial}{\partial A_X^j}$ ($j=1,\dots,6$) in place of $a_X$, $b_X$, $c_X$. From \eqref{6jss}, it is clear that for the given $f_i$, the $6j$-symbol can only be non-zero if the highest weights of the representations are as follows:

$$V^1=[1,1,0,0], V^2=[1,0,0,0], V^3=[1,1,0,0],$$  
$$ U=[1,0,0,0], W=[1,0,0,0], H=[1,0,0,0].$$ 

Thus, the $6j$-symbol is defined; let us compute its value.  

Note that we can reduce the sets of variables $Z^1, Z^2, Z^3, Z^4$ by keeping only those that arise in the decomposition of the given $f_i$.  

Now, let us describe the shifted lattice $H \cap pr^{-1}(D)$. Take formula \eqref{6jss} and consider the factors $f_1, \dots, f_4$ on the right-hand side. For convenience, as in the beginning of Section \ref{poror}, we replace the differential operators $\frac{\partial}{\partial A^j}$ in $f_i$ with the corresponding variables $A^j$ (unlike in Section \ref{poror}, we do not introduce an additional index $i$). The shifted lattice $H$ can be viewed as the set of exponents of monomials in the product $f_1 \cdots f_4$, where each factor is treated as a function of the variables $Z^1, Z^2, Z^3, Z^4$.  

A monomial in the product $f_1 \cdots f_4$ is a quadruple of monomials taken from $f_1, \dots, f_4$, respectively. When intersecting with $pr^{-1}(D)$, we only retain those quadruples that satisfy the following property: upon transitioning from the variables $Z^1, \dots, Z^4$ to $A^j_X$, the resulting quadruple of monomials must obey (cf. \eqref{6jss}):  

\begin{enumerate}  
\item $A^1_X$ appears with the same power in the first and fourth monomials.  
\item $A^2_X$ appears with the same power in the first and third monomials.  
\item $A^3_X$ appears with the same power in the second and third monomials.  
\item $A^4_X$ appears with the same power in the first and second monomials.  
\item $A^5_X$ appears with the same power in the second and fourth monomials.  
\item $A^6_X$ appears with the same power in the third and fourth monomials.  
\end{enumerate}  

It is easy to verify that from $f_1, \dots, f_4$, we must take quadruples of monomials\footnote{this is a quadruple of monomials in variables  $Z^1,...,Z^4$} of the form:  

\begin{equation}  
\label{aaaa}  
[A_{i,j}^1 A_{k}^2 A_{l}^4], \quad [A^4_{l}A^3_{i,j}A^5_{k}], \quad [A_{k}^2A_{i,j}^3A_l^6], \quad [A_{i,j}^1A^6_{l}A^5_{k}],  
\end{equation}  

where $(i,j,k,l)$ is a permutation of $1, \dots, 4$. There are $4!$ such quadruples. If the permutation $\sigma = (i,j,k,l)$ has sign $(-1)^{\sigma}$, then the listed monomials enter $f_1, \dots, f_4$ with coefficients:  

\begin{equation}  
\label{ssss}  
(-1)^{\sigma}, \quad -(-1)^{\sigma}, \quad (-1)^{\sigma}, \quad -(-1)^{\sigma}.  
\end{equation}  

Let us proceed to compute $\mathcal{J}_{\gamma}(z; L)$. According to the previous reasoning, the sum in \eqref{fj} runs over products of quadruples of monomials of the form \eqref{aaaa}. Now, let us find the coefficient for such a product.  

When treating $f_1, \dots, f_4$ as functions of the variables $Z^1, Z^2, Z^3, Z^4$, the monomials \eqref{aaaa} enter with coefficients $1$. After applying the projection $pr$, we also obtain monomials where the variables appear with power $1$. Thus, the coefficient $\frac{\Gamma(pr(x)+1)}{x!}$ for each product of monomials \eqref{aaaa} is $1$. Next, substituting \eqref{ssss} in place of the monomials \eqref{aaaa}, we obtain $1$. As a result, we get a sum of $1$ repeated $4!$ times.  

Thus, the $6j$-symbol in question equals $4!$.

\end{document}